\documentclass{amsart}
\usepackage{csquotes}
\usepackage{graphicx,amsaddr}
\usepackage{amsfonts}
\usepackage{amsmath}
\usepackage{mathtools}
\usepackage{amsthm}
\usepackage{amsfonts}
\usepackage{amssymb}
\usepackage{calc}
\usepackage{rotating}
\usepackage[all,cmtip]{xy}
\usepackage{pifont}
\usepackage{enumerate}
\usepackage{enumitem}
\usepackage{lscape}
\usepackage{array}

\theoremstyle{plain}
\newtheorem{defn}{Definition}[section]
\newtheorem{thm}[defn]{Theorem}
\newtheorem{cor}[defn]{Corollary}

\newtheorem{rem}[defn]{Remark}
\newtheorem{lem}[defn]{Lemma}

\numberwithin{equation}{section}

\title{On the convexity of spatial numetical range in normed algebras}

\author{H. V. Dedania}

\address{Dept. of Mathematics, Sardar Patel University, Vallabh Vidyanagar 388120, Gujarat, India}
\email{hvdedania@gmail.com}

\author{A. B. Patel*}

\address{Dept. of Mathematics, Sardar Patel University, Vallabh Vidyanagar 388120, Gujarat, India}
\email{avadhpatel663@gmail.com}
\newcolumntype{C}[1]{>{\centering\arraybackslash}m{#1}}

\begin{document}

\subjclass[2010]{Primary 46H05; Secondary \textbf{47A12}.}

\keywords{Banach algebra, Spatial Numerical Range.}

\begin{abstract}
In this article, we address the following question: Is it true that the spatial numerical range (SNR) $V_A(a)$ of an element $a$ in a normed algebra $(A, \|\cdot\|)$ is always convex? If the normed algebra is unital, then it is convex \cite[Theorem 3, P.16]{BoDu:71}. In non-unital case, we believe that the problem is still open and its answer seems to be negative. In search of such a normed algebra, we have proved that the SNR $V_A(a)$ is convex in several non-unital Banach algebras.
\end{abstract}
\maketitle

\section{Introduction}

\noindent
Let $(A, \|\cdot\|)$ be a normed algebra. Let $S(A)= \{x \in A: \|x\|=1\}$ be the unit sphere of $A$ and $D_{A}(x) = \{\varphi \in A^*: \|\varphi\|=1=\varphi(x)\}$ for each $x \in S(A)$. Let $V_{A}(a;x) = \{\varphi(ax): \varphi \in D_{A}(x)\}$ for $a \in A$. Then $V_{A}(a) = \cup\{V_{A}(a;x): x \in S(A)\}$ is the \emph{spatial numerical range (SNR)} and $\nu_{A}(a) = \sup \{|\lambda|: \lambda \in V_{A}(a)\}$ is the \emph{spatial numerical radius} of $a$ in $(A, \|\cdot\|)$. Some basic properties of $V_A(a)$ and its applications are given in~\cite{DePa:22(a),DePa:23}. The SNR plays an important role in studying the spectral extension property in non-unital Banach algebras and their unitization~\cite{DePa:22}. The SNR and NR (numerical range) of bounded operators have been studied in~\cite{Ba:62,Bo:69,GaHu:89,Ta:00}. In this paper, we study the SNR of elements of a normed algebra.  We have been searching a (necessarily non-unital) normed algebra $A$ and an element $a \in A$ such that the SNR $V_A(a)$ is not a convex set. We strongly believe that such an example must exist. Because there does exist a linear operator $T$ on the Banach space $(\mathbb C^2, \|\cdot\|_\infty)$ such that its SNR is not a convex set~\cite[P.98]{BoDu:71} and~\cite[P.357]{NiSc:64}. Unfortunately, this operator does not help us to find an element $a$ in $A$ such that $V_A(a)$ is not convex. Here we have calculated the SNR in some Banach algebras. In all these Banach algebras, it happens to be convex.

\section{Main Results}

\noindent
First we state the following known result. It asserts that if $A$ is a unital normed algebra, then  the SNR $V_A(a)$ is convex for each $a \in A$. We would like to include its simple and short proof.

\begin{lem}
  Let $(A, \|\cdot\|)$ be a normed algebra with an identity $1_A$ such that $\|1_A\| = 1$. Then
  \begin{enumerate}
  \item $V_A(a; x)$ is convex for each $a \in A$ and $x \in S(A)$;
  \item \cite[Lemma 2, P. 15]{BoDu:71} $V_A(a) = V_A(a; 1_A) \, (a \in A)$;
  \item \cite[Theorem 3, P. 16]{BoDu:71} $V_A(a)$ is convex for each $a \in A$.
  \end{enumerate}
\end{lem}

\begin{proof}
(1) Let $\lambda_1, \lambda_2 \in V_A(a; x)$. Then $\lambda_1 = \varphi_1(ax)$ and $\lambda_2 = \varphi_2(ax)$ for some $\varphi_1, \varphi_2 \in D_A(x)$. Take $r \in [0, 1]$ and $\varphi = r\varphi_1 + (1-r)\varphi_2$. Then $\varphi \in D_A(x)$ and $r\lambda_1 + (1-r)\lambda_2 = \varphi(ax) \in V_A(a; x)$. Thus $V_A(a; x)$ is convex.\\
(2) By definition, $V_A(a; 1_A) \subseteq V_A(a)$. Now let $\lambda \in V_A(a)$. Then $\lambda = \varphi(ax)$ for some $x \in S(A)$ and $\varphi \in D_A(x)$. Define $\psi(y) = \varphi(yx) \, (y \in A)$. Then $\psi \in D_A(1_A)$ and $\lambda = \varphi(ax) = \psi(a) = \psi(a1_A) \in V_A(a; 1_A)$. Thus $V_A(a) \subseteq V_A(a; 1_A)$.\\
(3) Since $1_A \in S(A)$, the convexity of $V_A(a)$ follows from (1) and (2) above.
\end{proof}

\begin{thm}
Let $(A, \|\cdot\|_A)$ and $(B, \|\cdot\|_B)$ be normed algebras. Consider the normed algebra $(A \times B, \|\cdot\|_1)$ with the co-ordinatewise operations and the norm $\|(a, b)\|_1 = \|a\|_A + \|b\|_B$. Then, $a \in A$ and $b \in B$,
$$V_{A \times B}(a, b) = \{rz + (1-r)w: r \in [0, 1], z \in V_A(a), \, w \in V_B(b)\}.$$
In particular, if $V_A(a)$ and $V_B(b)$ are convex, then $V_{A \times B}(a, b)$ is convex.
\end{thm}

\begin{proof}
Note that the dual of a normed algebra $(A \times B, \|\cdot\|_1)$ is isometrically linear isomorphic to $(A^* \times B^*, \|\cdot\|_{\infty})$ by the map $\varphi \longmapsto (\varphi_1, \varphi_2)$, where $\varphi(x, y) = \varphi_1(x) + \varphi_2(y)$ for $x \in A$ and $y \in B$.

Let $\lambda \in V_{A \times B}(a, b)$. Then $\lambda = \varphi((a, b)(x, y)) = \varphi(ax, by) = \varphi_1(ax) + \varphi_2(by)$ for some $(x, y) \in S(A \times B)$ and $\varphi = (\varphi_1, \varphi_2) \in S((A \times B)^*) = S(A^* \times B^*)$. So that
$$\|x\|_A + \|y\|_B = 1 \text{ and } \varphi_1(x) + \varphi_2(y) = \varphi(x, y) = 1 = \|\varphi\| = \max\{\|\varphi_1\|, \|\varphi_2\|\}.$$
If $\|x\|_A = 0$, then $\|y\|_B = 1 = \varphi_2(y) = \|\varphi_2\|$ and hence $\lambda = \varphi_2(by) \in V_B(b)$. Similarly, if $\|y\|_B = 0$, then $\lambda \in V_A(a)$. Now suppose that $\|x\|_A \neq 0 \neq \|y\|_B$. Take $\varphi_1(x) = r_1 + ir_2, \, x' = \frac{x}{\|x\|_A}, \, y' = \frac{y}{\|y\|_A}$, and $\|x\|_A = r \in (0, 1)$. Then
$$\varphi_2(y) = 1-r_1 - ir_2, \, \|y\|_B  = 1-r, \, |\varphi_1(x)| \leq r, \, \text{ and } \, |\varphi_2(y)| \leq 1-r.$$
If $r > r_1$, then $|1-r_1 -ir_2| = |\varphi_2(y)| \leq 1-r < 1-r_1$, which is not possible. If $r < r_1$, then $r < |r_1 + ir_2| =|\varphi_1(x)| \leq r$, which is also impossible. So $r = r_1$ and $r_2 = 0$. Hence $\varphi_1(x) = r_1 = r = \|x\|_A$ and $\varphi_2(y) = 1-r_1 = 1-r = \|y\|_B$. Now
$$\lambda = \varphi_1(ax) + \varphi_2(by)= r \varphi_1(ax') + (1-r)\varphi_2(by') \in rV_A(a; x') + (1-r)V_B(b; y').$$
Thus $V_{(A \times B)}(a, b) \subseteq \{rz + (1-r)w: r \in [0, 1], \, z \in V_A(a), \, w \in V_B(b)\}$.

Conversely, let $\lambda = rz + (1-r)w$ for some $r \in [0, 1], \, z \in V_A(a)$, and $w \in V_B(b)$. So, by definition, $z = \varphi_1(ax)$ and $w = \varphi_2(by)$ for some $x \in S(A), \, y \in S(B), \, \varphi_1 \in D_A(x)$, and $\varphi_2 \in D_B(y)$. Take $x' = rx, \, y' = (1-r)y$, and $\varphi = (\varphi_1, \varphi_2)$. Then $(x', y') \in S(A \times B)$, $\varphi \in D_{A \times B}(x', y')$, and
$$\lambda = rz + (1-r)w = \varphi_1(ax') + \varphi_2(by') = \varphi((a, b)(x', y')) \in V_{A \times B}(a, b).$$
Thus $\{rz + (1-r)w: r \in [0, 1], \, z \in V_A(a), \, w \in V_B(b)\} \subseteq V_{(A \times B)}(a, b)$.
\end{proof}

\noindent
We shall need the following result in the proof of Theorem \ref{p1}.

\begin{lem}\label{convex}\cite[Lemma 3]{TaUc:02}
Let $K \subset \mathbb C$ be bounded. Then the convex hull
$$co(K) = \{\sum_{i=1}^n\lambda_iz_i:z_i \in K, \, \lambda_i \geq 0, \, \sum_{i=1}^n \lambda_i = 1\}$$
and the $\sigma$-convex hull
$$\sigma co(K) = \{\sum_{i=1}^{\infty} \lambda_iz_i: z_i \in co(K), \, \lambda_i \geq 0, \, \sum_{i=1}^{\infty} \lambda_i =1\}$$
are identical, i.e., $co(K) = \sigma co(K)$.
\end{lem}

\noindent
Let $1 \leq p < \infty$, let $X$ be any non-empty set, and let $\omega: X \longrightarrow [1, \infty)$. Define $\ell^p(X, \omega) = \{f: X \longrightarrow \mathbb C : f\omega \in \ell^p(X)\}$ and $\|f\|_{p \omega} = (\sum_{x \in X} |f(x)|^p \omega(x)^p)^{\frac{1}{p}}$. Then the Banach space $(\ell^p(X, \omega), \|\cdot\|_{p\omega})$ is a semisimple, commutative Banach algebra with respect to the pointwise product; it is unital iff $X$ is finite.

\noindent
Let $X$ be a locally compact, Hausdorff space. Let $C_0(X)$ be the commutative $C^*$-algebra of all continuous functions vanishing at infinity. In \cite[Theorem 6.1]{GaKo:91}, Gaur and Kovarik proved that $V_{C_0(X)}(f) = co(f(X)) \, (f \in C_0(X))$. The following result is an analogue of this result.

\begin{thm}\label{p1}
Let $1 \leq p < \infty$, let $X$ be a non-empty set, let $\omega : X \longrightarrow [1, \infty)$, and let $A = (\ell^p(X, \omega), \|\cdot\|_{p\omega}, \cdot)$. Then
$$V_A(f) = co(f(X)) \quad{(f \in A)}.$$
In particular, $V_A(f)$ is convex.
\end{thm}

\begin{proof}
Suppose that $X$ is infinite. Let $f \in A$. So we can write $f=\sum_{j=1}^{\infty}f(x_j)\delta_{x_j}$, where $x_j \in X$. Let $\lambda \in co(f(X))$. Then $\lambda = \sum_{j=1}^{n}\lambda_j f(x_j)$ for some $n \in \mathbb N$, $\lambda_j \geq 0$, and $\sum_{j=1}^{n}\lambda_j = 1$. Let $g = \sum_{j=1}^{n}\frac{\lambda_j^{1/p}}{\omega(x_j)}\delta_{x_j} \in S(A)$ and $\phi_h \cong h= \sum_{j=1}^{n} \lambda_j^{1/q} \omega(x_j)\delta_{x_j} \in \ell^q(X, \frac{1}{\omega}) \cong A^*$, where $\frac{1}{p}+\frac{1}{q}=1$. Clearly $\|h\|_{q\frac{1}{\omega}}=1$ and $\phi_h(g) = \langle g, h \rangle = \sum_{j=1}^{n} g(x_j)h(x_j) = \sum_{j=1}^{n} \lambda_j = 1$.
Thus $\phi_h \in D_A(g)$. Hence
\begin{eqnarray*}
\lambda & = &\sum_{j=1}^{n} \lambda_j f(x_j) = \langle \sum_{j=1}^{n} \frac{{\lambda_j}^{1/p}f(x_j)}{\omega(x_j)}\delta_{x_j}, \sum_{j=1}^{n} \lambda_j^{1/q}\omega(x_j) \delta_{x_j}\rangle \\
& = & \langle \sum_{j=1}^{\infty} f(x_j)\delta_{x_j} \cdot \sum_{j=1}^{n} \frac{{\lambda_j}^{1/p}}{\omega(x_j)}\delta_{x_j}, \sum_{j=1}^{n} \lambda_j^{1/q}\omega(x_j) \delta_{x_j}\rangle\\
& = & \langle f \cdot g, h \rangle = \phi_h(f \cdot g) \in V_A(f;g) \subseteq V_A(f).
\end{eqnarray*}
Thus $co(f(X)) \subseteq V_A(f)$.\\
Conversely, let $\lambda \in V_A(f)$. Then there exists $g \in S(A)$ such that $\lambda \in V_A(f;g)$. So there exists $\varphi \in D_A(g)$ such that $\lambda = \varphi(f \cdot g)$. Since $\varphi \in D_A(g) \subseteq A^* \cong \ell^q(X, \frac{1}{\omega})$, where $\frac{1}{p}+\frac{1}{q} = 1$, there exists $h \in \ell^q(X, \frac{1}{\omega})$ such that $\|h\|_{q \frac{1}{\omega}}=1$ and $\varphi = \phi_h$. Let $g=\sum_{j=1}^{\infty}g(x_j)\delta_{x_j}$ and $h \sim \sum_{j=1}^{\infty}h(x_j)\delta_{x_j}$. Let $g(x_j)h(x_j)=\alpha_j+i\beta_j$, where $\alpha_j, \beta_j \in \mathbb R$.  Then
\begin{eqnarray*}
1 & = & \varphi(g) = \phi_h(g) = |\phi_h(g)|
 =  |\sum_{j=1}^{\infty} g(x_j)h(x_j)|\\
& \leq & \sum_{j=1}^{\infty} |g(x_j)h(x_j)| =  \sum_{j=1}^{\infty} |g(x_j) \omega(x_j) h(x_j)\frac{1}{\omega(x_j)}|\\
& \leq & \|g \omega\|_p \|h \frac{1}{\omega}\|_q = \|g\|_{p \omega} \|h\|_{q \frac{1}{\omega}} =1.
\end{eqnarray*}
Therefore
\begin{eqnarray}\label{p1e1}
1=\sum_{j=1}^{\infty}|g(x_j)h(x_j)| = \sum_{j=1}^{\infty}|\alpha_j+i\beta_j|.
\end{eqnarray}
Also we have
\begin{eqnarray}\label{p1e2}
1 = \phi_h(g) = \sum_{j=1}^{\infty}g(x_j)h(x_j) = \sum_{j=1}^{\infty}\alpha_j+i\sum_{j=1}^{\infty}\beta_j.
\end{eqnarray}
Therefore, by equations~(\ref{p1e1}) and (\ref{p1e2}), we have
$$\sum_{j=1}^{\infty}|\alpha_j+i\beta_j| = \sum_{j=1}^{\infty}  \alpha_j$$
Therefore $\beta_j=0$, for all $j$. i.e, $g(x_j)h(x_j) \in \mathbb R$. But then $\sum_{j=1}^{\infty}|g(x_j)h(x_j)|=1$, we must have $g(x_j)h(x_j) \in [-1, 1]$, for all $j$.
If $g(x_j)h(x_j) < 0 $ for some j, then we get $1=\sum_{j=1}^{\infty}g(x_j)h(x_j) < \sum_{j=1}^{\infty}|g(x_j)h(x_j)|=1$, which is a contradiction. Therefore we have $g(x_j)h(x_j) \in [0,1]$ for every $j$. Now take $\lambda_j = g(x_j)h(x_j)$ ($j \in \mathbb N$). Then each $\lambda_j \geq 0, \, \sum_{j=1}^{\infty}\lambda_j=1$, and
\begin{eqnarray*}
\lambda & = & \varphi(f \cdot g) = \phi_h(f\cdot g)=\langle f \cdot g, h \rangle\\
& = & \sum_{j=1}^{\infty} f(x_j)g(x_j)h(x_j)
 =  \sum_{j=1}^{\infty} \lambda_j f(x_j) \in \sigma co(f(X)).
\end{eqnarray*}
Thus $V_A(f) \subseteq \sigma co(f(X))$. Therefore, by Lemma \ref{convex}, we get $V_A(f) \subseteq co(f(X))$.\\
Suppose $X$ is finite. Then the proof is much easier and similar to the above proof.
\end{proof}

\begin{cor}
Let $1 \leq p < r < \infty$ and $A = (\ell^p(X, \omega), \|\cdot\|_{r\omega}, \cdot)$. Then $V_A(f)$ is convex for each $f \in A$.
\end{cor}

\begin{proof}
Let $B = (\ell^r(X, \omega), \|\cdot\|_{r\omega})$. Since $A$ is dense in $B$, $V_A(f) \subseteq V_B(f)$. By Theorem~\ref{p1}, $V_B(f) = co(f(X))$. As per the arguments in the proof of Theorem~\ref{p1}, we can show that $co(f(X)) \subseteq V_A(f)$. Hence $V_A(f) = co(f(X))$ is convex.
\end{proof}

\noindent
\textbf{Notation:} Throughout, $\mathbb D = \{z \in \mathbb C: |z| \leq 1\}$.\\

\noindent
Let $(S, \cdot)$ be a semigroup. Then a \emph{weight} on $S$ is a strictly positive function $\omega$ on $S$ such that $\omega(s \cdot t) \leq \omega(s)\omega(t) \, (s, t \in S)$. Now we calculate the SNR $V_{\ell^1(S, \omega)}(f)$ in the weighted discrete semigroup algebra $(\ell^1(S, \omega), \|\cdot\|_{1\omega})$ with the \emph{convolution product $\star$}; see~\cite[P.159]{Da:00} for its detailed definition. In next result, we show that, under some mild conditions on $S$ and $\omega$, the SNR $V_{\ell^1(S, \omega)}(f)$ is a closed convex.

\begin{thm}\label{S}
Let $S$ be a right cancellative semigroup without right identity. Let $\omega : S \longrightarrow [1, \infty)$ be a weight on $S$ such that $\omega(t) = 1$  for some $t \in S$. Then
$$V_{\ell^1(S, \omega)}(f) = \|f\|_{1\omega}\mathbb D \quad{(f \in \ell^1(S, \omega))}.$$
\end{thm}

\begin{proof}
Let $A = \ell^1(S, \omega)$ and $f \in A$. Let $f = \sum_{n=1}^{\infty} f(s_n)\delta_{s_n}$. Fix $t \in S$ such that $\omega(t) = 1$. Since $S$ is right cancellative and without right identity, $t \neq s \cdot t \; (s \in S)$. Let $|z| \leq \|f\|_{1\omega}$. Take $g = \delta_t \, \in S(A)$ and $\phi_h \sim h = \delta_t + \sum_{n=1}^{\infty} h(s_n) \delta_{s_n \cdot t} \in \ell^{\infty}(S, \frac{1}{\omega}) \cong A^*$, where $h(s)=
\frac{z \overline{f(s)}\omega(s)}{\|f\|_{1\omega} |f(s)|} \mbox{if } s \in suppf$ and $h(s) = 0$ otherwise.
Then $h(t) = 1$, $|h(s)| \leq \omega(s) \, (s \in S)$, and hence $\|h\|_{\infty}\frac{1}{\omega} = 1$. Also we have
$\phi_h(g) = \langle  g, h \rangle =  \langle \delta_t , \delta_t + \sum_{n=1}^{\infty} h(s_n) \delta_{s_n \cdot t} \rangle =1$. Now
\begin{eqnarray*}
\phi_h(f \star g) & = & \langle f \star g, h \rangle
 =  \langle f \star \delta_t, h \rangle \\
& = & \langle \sum_{n=1}^{\infty} f(s_n)\delta_{s_n \cdot t},\delta_t + \sum_{n=1}^{\infty} h(s_n)\delta_{s_n \cdot t} \rangle =  \sum_{n=1}^{\infty}f(s_n)h(s_n)\\
 & = & \sum_{n=1}^{\infty} \frac{z}{\|f\|_{1 \omega}}|f(s_n)|\omega(s_n) = z.
\end{eqnarray*}
Thus $z = \phi_h(f \star g) \in V_A(f;g) \subseteq V_A(f)$. On the other hand, if $z \in V_A(f)$, then $|z| \leq \|f\|_{1\omega}$ is always true. This completes the proof.		
\end{proof}

\begin{rem}
We do not know whether the condition ``$\omega(t) = 1$  for some $t \in S$'' in above theorem is required or not? But the conditions on the semigroup $S$ can not be omitted as the following examples suggest:\\
(i) Consider the semigroup $(\mathbb Z^+, +)$, which is unital. Then $V_{\ell^1(\mathbb Z^+)}(\delta_0) = \{1\}$.\\
(ii) Let $\mathbb N_r = \mathbb N$ with the binary operation $m \cdot n = n$. Then $\mathbb N_r$ is not right cancellative and $V_{\ell^1(\mathbb N_r)}(\delta_n) = \{1\}$ for any $n \in \mathbb N_r$.
\end{rem}

\begin{defn}
Let $G$ be a dense subgroup of $\mathbb R$. Let $G_1 = G \cap (0, 1)$. For $f, g \in \ell^1(G_1)$, define $f \star g: G_1 \longrightarrow \mathbb C$ as
\begin{eqnarray*}
(f \star g)(s)&=& \sum \{f(u)g(u): u, v \in G_1 \text{ and } s = u+v\} \\
   &=& \sum \{f(s-t)g(t): t \in G_1 \text{ and } t < s\}.
\end{eqnarray*}
Then \index{$\ell^1(G_1)$}$(\ell^1(G_1), \star, \|\cdot\|_1)$ is a commutative radical Banach algebra. It is called the \emph{discrete volterra algebra} on $G_1$.
\end{defn}

\noindent
Next we show that the SNR is convex in $\ell^1(G_1)$.

\begin{thm}\label{thm 3.3.1}
Let $f \in A = \ell^1(G_1)$. Then $V_{A}(f) = \|f\|_1 \mathbb D$; in particular, $V_A(f)$ is convex.
\end{thm}

\begin{proof}
Let $f \in A$. Then the $supp f$ is countable. So we can write $f = \sum_{n=1}^{\infty}f(s_n)\delta_{s_n}$, where $s_n \in G_1$. Let $z \in \|f\|_1\mathbb D^o$ i.e., $|z| < \|f\|_1$. Since $\|f\|_1 = \sum_{n=1}^{\infty}|f(s_n)|$, there exists $n_0 \in \mathbb N$ such that $|z| < \sum_{n = 1}^{n_0}|f(s_n)| = \|\widetilde{f}\|_1$(say). Choose $t \in G_1$ such that $s_n + t < 1$. Now consider $g = \delta_t$. Let $h = \delta_t + \sum_{n=1}^{n_0}\frac{z\overline{f(s_n)}}{\|\widetilde{f}\|_1|f(s_n)|}\delta_{s_n+t}$. Then clearly, $\|g\|_1 = 1 = \|h\|_{\infty}$ and $\phi_h(g) = 1$. Now
\begin{eqnarray*}
  \phi_h(f \star g) &=& \langle f \star g, h \rangle = \sum_{s \in G_1}(f \star g)(s)h(s) \\
   &=& \sum_{s \in G_1}\bigl(\sum_{\underset{u < s}{u \in G_1}}f   (s-u)g(u)\bigr)h(s)\\
   &=& \sum_{\underset{t<s}{s \in G_1}}f(s-t)h(s) \quad{(\because g = \delta_t)} \\
   &=& \sum_{\underset{s_n+t<1}{n \leq n_0}}f(s_n)h(s_n+t) \quad{(\because s-t = s_n)}\\
   &=& \sum_{n = 1}^{n_0}\frac{z|f(s_n)|}{\|\widetilde{f}\|_1}
   = \frac{z}{\|\widetilde{f}\|_1}\sum_{n = 1}^{n_0}|f(s_n)|
  = z.
\end{eqnarray*}
Therefore $z \in V_A(f;g)$ and so, $\|f\|_1 \mathbb D^o \subseteq V_{A}(f;g)$. By \cite[Theorem 2.1(2)]{DePa:23}, $V_{A}(f;g)$ is compact. Therefore we get the inclusion $\|f\|_1\mathbb D \subseteq V_A(f;g)$. By definition, $V_A(f;g) \subseteq V_A(f)$ and so $\|f\|_1\mathbb D \subseteq V_A(f)$. Since $\nu_A(a) \leq \|a\|$ is always true in any normed algebra, we get $V_A(f) \subseteq \|f\|_1 \mathbb D$. Therefore $V_{A}(f) = \|f\|_1 \mathbb D$. Hence $V_{A}(f)$ is convex.
\end{proof}

\noindent
Now we calculate the SNR for the elements of $L^1(\mathbb R)$ and $L^1[0, 1]$. We could calculate it only for non-negative functions whose support does not contain zero.
\begin{lem}
Let $A = L^1(\mathbb{R})$. Let $f \in A$ such that $f \geq 0$ and $0 \notin suppf$. Then $V_A(f) = \|f\|_1\mathbb D$. In particular, $V_{A}(f)$ is convex.
\end{lem}

\begin{proof}
Let $f \in A$. Since $\nu_A(f) \leq \|f\|$ is always true, we get $V_A(f) \subseteq \|f\|_1 \mathbb D$. For the reverse inclusion, let $z \in \mathbb{D}$. Since $0 \notin suppf$ and $supp f$ is a closed set, there exists $a > 0$ such that $supp f \subseteq \mathbb{R} \setminus (-a, a)$.
Let $g = \frac{2}{a}\chi_{_{[0, \frac{a}{2}]}} \in S(A)$ and $h = \chi_{_{[0, \frac{a}{2}]}} + z \chi_{_{\mathbb{R} \setminus [\frac{-a}{2}, \frac{a}{2}]}} \in L^{\infty}(\mathbb R)$. We know that $L^1(\mathbb{R})^* = L^{\infty}(\mathbb{R})$. So that $\phi_h \cong h \in A^*$ and $\|h\|_{\infty} = 1$. Clearly, $\phi_h(g) = \int_{\mathbb R} g(t)h(t) dt = 1$ and so $\varphi_h(f \star g) \in V_A(f)$. Now it is enough to show that $\phi_h(f \star g) = \|f\|_1z$. In fact,
\begin{eqnarray*}
 \phi_h(f \star g) &=& \int_{\mathbb R} \biggl[ \int_{\mathbb R} f(s-t)g(t) dt \biggr] h(s) ds\\
   &=& \int_{\mathbb R} \biggl[\int_{\mathbb R} f(s-t)g(t) dt\biggr] \chi_{_{[0, \frac{a}{2}]}}(s) ds\\
   & & + \int_{\mathbb R} \biggl[\int_{\mathbb R} f(s-t)g(t) dt\biggr] z \chi_{_{\mathbb{R} \setminus [\frac{-a}{2}, \frac{a}{2}]}}(s) ds\\
    &=&  \frac{2}{a} \int_{\mathbb R} \biggl[\int_{0}^{\frac{a}{2}} f(s-t) dt\biggr]\chi_{_{[0, \frac{a}{2}]}}(s) ds\\
   & & + z\int_{\mathbb R} \biggl[\int_{\mathbb R} f(s-t)\chi_{_{\mathbb{R} \setminus [\frac{-a}{2}, \frac{a}{2}]}}(s)ds\biggr] g(t) dt
\end{eqnarray*}

\begin{eqnarray*}
   &=& \frac{2}{a} \int_{0}^{\frac{a}{2}} \biggl[\int_{0}^{\frac{a}{2}} f(s-t) dt\biggr] ds + \frac{2z}{a}\int_{0}^{\frac{a}{2}} \biggl[\int_{-\infty}^{\frac{-a}{2}} f(s-t)ds + \int_{\frac{a}{2}}^{\infty} f(s-t)ds\biggr] dt\\
  % &=& \frac{2}{a} \int_{0}^{\frac{a}{2}} \biggl[\int_{s-\frac{a}{2}}^{s} f(x) dx\biggr] ds\\
%     & &+  \frac{2z}{a} \int_{0}^{\frac{a}{2}} \biggl[\int_{-\infty}^{-\frac{a}{2}-t} f(x) dx + \int_{\frac{a}{2}-t}^{\infty} f(x) dx\biggr] dt\\
      &=& 0 +  \frac{2z}{a} \int_{0}^{\frac{a}{2}} \biggl[\int_{-\infty}^{-\frac{a}{2}-t} f(x) dx + \int_{\frac{a}{2}-t}^{\infty} f(x) dx\biggr] dt \quad{(\because suppf \subseteq \mathbb{R} \setminus (-a, a))}\\
   &=& \frac{2z}{a}\|f\|_1 \int_{0}^{\frac{a}{2}} dt = z\|f\|_1.
\end{eqnarray*}
Hence $V_A(f) = \|f\|_1\mathbb D$.
\end{proof}

\begin{defn}
Define $L^1[0, 1] = \{f: [0, 1] \longrightarrow \mathbb C: \int_{0}^{1}|f(t)|dt< \infty\}$. For $f, g \in L^1[0, 1]$, define
$$(f \star g)(s) = \int_{0}^{s}f(s-t)g(t)dt \text{ and }
  \|f\|_1 = \int_{0}^{1}|f(t)| dt.$$
Then $(L^1[0, 1], \star, \|\cdot\|_1)$ is a Banach algebra. It is called the \emph{Volterra algebra}.
\end{defn}

\begin{lem}
Let $A = L^1[0, 1]$. Let $f \in A$ such that $f \geq 0$ and $0, 1 \notin suppf$. Then $V_A(f) = \|f\|_1\mathbb D$. In particular, $V_{A}(f)$ is convex.
\end{lem}

\begin{proof}
Since $0, 1 \notin suppf$, there exists $\delta > 0$ such that $suppf \subseteq [\delta, 1-\delta]$. Because $\nu_A(a) \leq \|a\|$ is always true, we get $V_A(f) \subseteq \|f\|_1 \mathbb D$. For the reverse inclusion, let $z \in \mathbb{D}$. Let $g = \frac{2}{\delta}\chi_{_{[0, \frac{\delta}{2}]}} \in S(A)$ and $h = \chi_{_{[0, \frac{\delta}{2})}} + z \chi_{_{[\frac{\delta}{2}, 1]}}$. We know that $L^1([0, 1])^* = L^{\infty}([0, 1])$ so $\phi_h \cong h \in A^*$ and $\|h\|_{\infty} = 1$. Clearly, $\phi_h(g) = \int_{0}^{1} g(t)h(t) dt = 1$ and so $\varphi_h(f \star g) \in V_A(f)$. Now it is enough to show that $\phi_h(f \star g) = \|f\|_1z$.
\begin{eqnarray*}
\phi_h(f \star g) &=& \int_{0}^{1} \biggl[ \int_{0}^{s} f(s-t)g(t) dt \biggr] h(s) ds\\
  & = & \frac{2}{\delta}\int_{0}^{1} \biggl[\int_{0}^{s} f(s-t)\chi_{_{[0, \frac{\delta}{2})}} dt\biggr] \chi_{_{[0, \frac{\delta}{2}]}}(s) ds\\
  & & + \frac{2z}{\delta} \int_{0}^{1} \biggl[\int_{0}^{s} f(s-t)\chi_{_{[0, \frac{\delta}{2})}} dt\biggr] \chi_{_{[\frac{\delta}{2}, 1]}}(s) ds\\
    %&=& \frac{2}{\delta}\int_{0}^{\frac{\delta}{2}} \biggl[\int_{0}^{s} f(s-t)\chi_{_{[0, \frac{\delta}{2})}} dt\biggr] ds\\
%   & & +\frac{2z}{\delta} \int_{0}^{1} \biggl[\int_{0}^{s} f(s-t)\chi_{_{[0, \frac{\delta}{2})}} dt\biggr] \chi_{_{[\frac{\delta}{2}, 1]}}(s) ds\\
   &=&  0 + \frac{2z}{\delta} \int_{0}^{1} \biggl[\int_{0}^{s} f(s-t)\chi_{_{[0, \frac{\delta}{2})}}(t) dt\biggr] \chi_{_{[\frac{\delta}{2}, 1]}}(s) ds \, (\because suppf \subseteq [\delta, 1-\delta])\\
   %&=& \frac{2z}{\delta} \int_{\frac{\delta}{2}}^{1} \biggl[\int_{0}^{s} f(s-t)\chi_{_{[0, \frac{\delta}{2})}}(t) dt\biggr] ds\\
   &=& \frac{2z}{\delta} \int_{\frac{\delta}{2}}^{1} \biggl[\int_{0}^{\frac{\delta}{2}} f(s-t) dt\biggr] ds
   = \frac{2z}{\delta} \int_{0}^{\frac{\delta}{2}} \biggl[ \int_{\frac{\delta}{2}}^{1} f(s-t) ds\biggr] dt\\
   &=& \frac{2z}{\delta} \int_{0}^{\frac{\delta}{2}} \biggl[ \int_{\frac{\delta}{2}-t}^{1-t} f(x) dx\biggr] dt \quad{(\because s-t = x)}\\
   &=& \frac{2z}{\delta} \int_{0}^{\frac{\delta}{2}} \|f\|_1 dt = z\|f\|_1.
\end{eqnarray*}
Hence $V_A(f) = \|f\|_1\mathbb D$.
\end{proof}

\noindent
Finally, we calculate the SNR in the Banach algebra $\mathbb C^2$ with 35 different products and algebra norms. In all these Banach algebras, the SNR happens to be convex. In fact, we considered more than 150 binary operations on $\mathbb C^2$. However, $\mathbb C^2$ becomes an algebra with only 35 binary operations. The calculation of $V_A(a)$ is standard.

\begin{table}[ht]
%\hspace{5cm}
The SNR of an element of $\mathbb C^2$ with different products and suitable norms\\ (see the ``note" below for the notations used in the table)
\small\addtolength{\tabcolsep}{-5pt}
\begin{tabular}{| c | C{4cm} | C{1.7cm} | C{6cm} |}
\hline
No. &  Product        & Norm & $V_A(a)$ \\[2ex]
\hline
1   &$(x_1y_1, 0)$  &  $\|\cdot\|_{p}$     & $a_1I$ \\
\hline
2   &$(x_2y_2, 0)$ &  $\|\cdot\|_{p}$ & $\{a_2r^{\frac{1}{q}}(1-r)^{\frac{1}{p}}e^{i\theta}: r \in I, \theta \in I_{\pi}\}$ \\
\hline
3   &$(0, x_1y_1)$  &$\|\cdot\|_{p}$ & $\{a_1r^{\frac{1}{p}}(1-r)^{\frac{1}{q}}e^{i\theta}: r \in I, \theta \in I_{\pi}\}$  \\
\hline
4   &$(0, x_2y_2)$  &$\|\cdot\|_{p}$  & $a_2I$ \\
\hline
5   &$(2x_1y_1, 0)$    & $\|\cdot\|_{p2}$ & $2a_1I$ \\
\hline
6   &$(2x_2y_2, 0)$  & $2\|\cdot\|_{p}$ & $\{2a_2r^{\frac{1}{q}}(1-r)^{\frac{1}{p}}e^{i\theta}: r \in I, \theta \in I_{\pi}\}$ \\
\hline
7   &$(0, 2x_1y_1)$  & $2\|\cdot\|_{p}$ & $\{2a_1r^{\frac{1}{p}}(1-r)^{\frac{1}{q}}e^{i\theta}: r \in I, \theta \in I_{\pi}\}$   \\
\hline
8   &$(0, 2x_2y_2)$  & $2\|\cdot\|_{p}$ & $2a_2I$ \\
\hline
9   &$(3x_1y_1, 0)$   & $3\|\cdot\|_{p}$ & $3a_1I$ \\
\hline
10   &$(3x_2y_2, 0)$  & $3\|\cdot\|_{p}$  & $\{3a_2r^{\frac{1}{q}}(1-r)^{\frac{1}{p}}e^{i\theta}: r \in I, \theta \in I_{\pi}\}$ \\
\hline
11  &$(0, 3x_1y_1)$   & $3\|\cdot\|_{p}$  &$\{3a_1r^{\frac{1}{p}}(1-r)^{\frac{1}{q}}e^{i\theta}: r \in I, \theta \in I_{\pi}\}$   \\
\hline
12   &$(0, 3x_2y_2)$ & $3\|\cdot\|_{p}$  & $3a_2I$ \\
\hline
13   &$(x_1y_1, x_1y_1)$  &$2^{\frac{1}{p}}\|\cdot\|_{p}$ & $\{a_1r + a_1r^{\frac{1}{p}}(1-r)^{\frac{1}{q}}e^{i\theta}: r \in I, \theta \in I_{\pi}\}$  \\
\hline
14   &$(x_1y_1, x_1y_2)$  & Any &  $\{a_1\}$  \\
\hline
15   &$(x_1y_1, x_2y_1)$   & $\|\cdot\|_{p}$ & $\{a_1r + a_2r^{\frac{1}{p}}(1-r)^{\frac{1}{q}}e^{i\theta}: r \in I, \theta \in I_{\pi}\}$\\
\hline
16  &$(x_1y_1, x_2y_2)$ & $\|\cdot\|_{p}$ & $co\{a_1, a_2\}$  \\
\hline
17  &$(x_1y_2, x_2y_2)$  &$\|\cdot\|_{p}$ & $\{a_1r^{\frac{1}{q}}(1-r)^{\frac{1}{p}}e^{i\theta} + a_2(1-r): r \in I, \theta \in I_{\pi}\}$   \\
\hline
18   &$(x_2y_1, x_2y_2)$   & Any & $\{a_2\}$  \\
\hline
19   &$(x_2y_2, x_2y_2)$  &$2^{\frac{1}{p}}\|\cdot\|_{p}$ & $\{a_2r^{\frac{1}{q}}(1-r)^{\frac{1}{p}}e^{i\theta} + a_2(1-r): r \in I, \theta \in I_{\pi}\}$   \\
\hline
20   &$(2x_1y_1, x_1y_1)$  &$3\|\cdot\|_{p}$ & $\{2a_1r + a_1r^{\frac{1}{p}}(1-r)^{\frac{1}{q}}e^{i\theta}: r \in I, \theta \in I_{\pi}\}$  \\
\hline
21   &$(2x_2y_2, x_2y_2)$  & $3\|\cdot\|_{p}$ & $\{2a_2r^{\frac{1}{q}}(1-r)^{\frac{1}{p}}e^{i\theta} + a_2(1-r): r \in I, \theta \in I_{\pi}\}$  \\
\hline
22   &$(x_1y_1, 2x_1y_1)$  &$3\|\cdot\|_{p}$ & $\{a_1r + 2a_1r^{\frac{1}{p}}(1-r)^{\frac{1}{q}}e^{i\theta}: r \in I, \theta \in I_{\pi}\}$  \\
\hline
23   &$(x_1y_1, 2x_2y_2)$ &$2\|\cdot\|_{p}$  & $\{a_1r + 2a_2(1-r): r \in I\}$  \\
\hline
24   &$(x_2y_2, 2x_2y_2)$  & $3\|\cdot\|_{p}$ & $\{a_2r^{\frac{1}{q}}(1-r)^{\frac{1}{p}}e^{i\theta} + 2a_2(1-r): r \in I, \theta \in I_{\pi}\}$  \\
\hline
25   &$(x_1y_2 + x_2y_1, x_2y_2)$  &$\|\cdot\|_{1}$ &  $a_2+a_1 \mathbb D$ \\
\hline
26   &$(x_1y_1, x_1y_2 + x_2y_1)$   &$\|\cdot\|_{1}$  & $a_1+a_2 \mathbb D$   \\
\hline
27   &$(x_1y_1+ x_1y_2 + x_2y_1 +x_2y_2, 0)$  & $\|\cdot\|_{1}$  & $|a_1+a_2| \mathbb D$  \\
\hline
28   &$(0, x_1y_1+ x_1y_2 + x_2y_1 +x_2y_2)$ & $\|\cdot\|_{1}$  & $|a_1+a_2| \mathbb D$  \\
\hline
29   &$(x_1y_1+ x_1y_2 + x_2y_1, x_2y_2)$  & $\|\cdot\|_{1}$ & $a_2+a_1 \mathbb D$   \\
\hline
30   &$(x_1y_2 + x_2y_1 + x_2y_2, x_2y_2)$  &  $2\|\cdot\|_{1}$ & $a_2 + (a_1 + a_2)\mathbb D$  \\
\hline
31   &$(x_1y_1 + x_2y_2, x_1y_2 +x_2y_1)$  & $\|\cdot\|_1$ &  $a_1+a_2 \mathbb D$  \\
\hline
32   &$(x_1y_1 + x_1y_2, x_2y_1 + x_2y_2)$ & $\|\cdot\|_1$ & $\{r(a_1+a_2e^{i\theta})+(1-r)(a_1e^{i\theta}+a_2): r \in I, \theta \in I_{\pi}\}$
\\
\hline
33   &$(x_1y_1-x_2y_2, x_1y_2+x_2y_1)$  & $\|\cdot\|_1$ & $a_1+a_2 \mathbb D$     \\
\hline
34   &$(x_1y_1 + x_2y_1, x_1y_2+x_2y_2)$  & $\|\cdot\|_1$ & $\{a_1+a_2\}$ \\
\hline
35   &$(x_1y_2 + x_2y_1, x_1y_1+x_2y_2)$  & $\|\cdot\|_1$ & $a_2+a_1 \mathbb D$  \\
\hline
\end{tabular}
\end{table}

\newpage
\noindent
\textbf{Note:} Following notations are used in the above table.
(1) $x = (x_1, x_2)$ and $y = (y_1, y_2)$ are elements of $\mathbb C^2$.
(2) $I = [0, 1]$ and $I_{\pi} = [-\pi, \pi)$.
(3) $1 \leq p, q \leq \infty$ and $\frac{1}{p} + \frac{1}{q} = 1$.
(4) $\|x\|_p = (|x_1|^p + |x_2|^p)^{\frac{1}{p}}$ and $\|x\|_{\infty} = \max\{|x_1|, |x_2|\}$.
(5) $c\|\cdot\|_p$ means $\|x\| = c\|x\|_p$.
(6) The SNR $V_A(a)$ is calculated for $a=(a_1, a_2) \in \mathbb C^2$.
%\begin{enumerate}
%\item $x = (x_1, x_2)$ and $y = (y_1, y_2)$ are elements of $\mathbb C^2$.
%\item $I = [0, 1]$ and $I_{\pi} = [-\pi, \pi)$.
%\item $1 \leq p, q \leq \infty$ and $\frac{1}{p} + \frac{1}{q} = 1$.
%\item $\|x\|_p = (|x_1|^p + |x_2|^p)^{\frac{1}{p}}$ and $\|x\|_{\infty} = \max\{|x_1|, |x_2|\}$.
%\item $c\|\cdot\|_p$ means $\|x\| = c\|x\|_p$.
%\item The SNR $V_A(a)$ is calculated for $a=(a_1, a_2) \in \mathbb C^2$.
%\end{enumerate}

\noindent
\textbf{Conflict of interest :} {On behalf of all authors, the corresponding author states that there is no conflict of interest.}

\noindent
\textbf{Acknowledgement :} {The second author is thankful to the Council of Scientific and Industrial Research (CSIR), New Delhi, for providing Senior Research Fellowship.}

\end{document}